%% file: paper.tex
\documentclass[10pt]{article}
\usepackage[utf8]{inputenc}
\usepackage[a4paper, width = 160mm, left = 25mm, top = 30mm, height = 240mm]{geometry}
\usepackage{amsmath}
\usepackage{amsfonts}
\usepackage{amssymb}
\usepackage{cite}
\usepackage{hyperref}
\usepackage{caption}
\usepackage{float}
\usepackage{colortbl}
\usepackage{graphicx}
\usepackage{mathtools}
\usepackage{empheq}
\usepackage{setspace}
\usepackage{widetext}
\usepackage{cases}
\usepackage{subfigure}

\title{\textbf{THE TAKAGI CURVE AND THE $\beta$-CANTOR FUNCTION FROM MECHANICAL LAWS}}
\author{Javier Rodríguez-Cuadrado* and Jesús San Martín\\
\normalsize Universidad Politécnica de Madrid, Madrid 28012, Spain\\
\normalsize \href{mailto:javier.rodriguez.cuadrado@upm.es}{javier.rodriguez.cuadrado@upm.es}, \href{mailto:jesus.sanmartin@upm.es}{jesus.sanmartin@upm.es}
}
\date{\today}

\begin{document}

\maketitle

\begin{abstract}
This work shows that fractals can be obtained from Mechanical Laws without being forced by any algorithm, closing the gap between the Platonic world of Mathematics and Nature. Fractal tree crown directly emerges when applying elasticity theory to branching stresses in a binary tree. Vertical displacements of nodes are given by the Takagi curve, while the horizontal ones are given by a linear combination of inverses of $\beta$-Cantor functions. In addition, both fractal dimensions are related, which suggests a deeper connection between the Takagi Curve and the $\beta$-Cantor function.
\\ \\
*Corresponding author
\end{abstract}

\vspace{0.1 cm}

\input{introduccion}
\input{modelo}
\input{vertical}
\input{horizontal}
\input{iteraciones}
\input{dimension}
\input{conclusiones}

\section*{Acknowledgements}

This action is financed by Universidad Politécnica de Madrid as part of the UPM-Funded Research, Development and Innovation Programme, specifically targeting Funding for predoctoral contracts for the completion of doctoral degrees at UPM schools, faculties, and R\&D centres and institutes including a period of at least three months as a visiting researcher outside Spain (international doctoral mention).

\end{document}

%% file: introduccion.tex
\section{INTRODUCTION}

Fractals are everywhere, so a fundamental question arises: what is the reason for their ubiquity? Significant development has been done from pure Mathematics \cite{1}, but where is the bridge between the real world and the Platonic world of Mathematics? Where is the link between the ubiquitous fractals we see in Nature and Mathematics? The answer to this first question is easy. As fractal structures are naturally occurring, they are obviously ruled by the laws of science. So our original question evolves towards this specific question: Why do fractals exist? or more specifically: what are exactly the real mechanisms underlying in its physical formation? 

\begin{figure}[h]
\centering
\includegraphics[width=8cm]{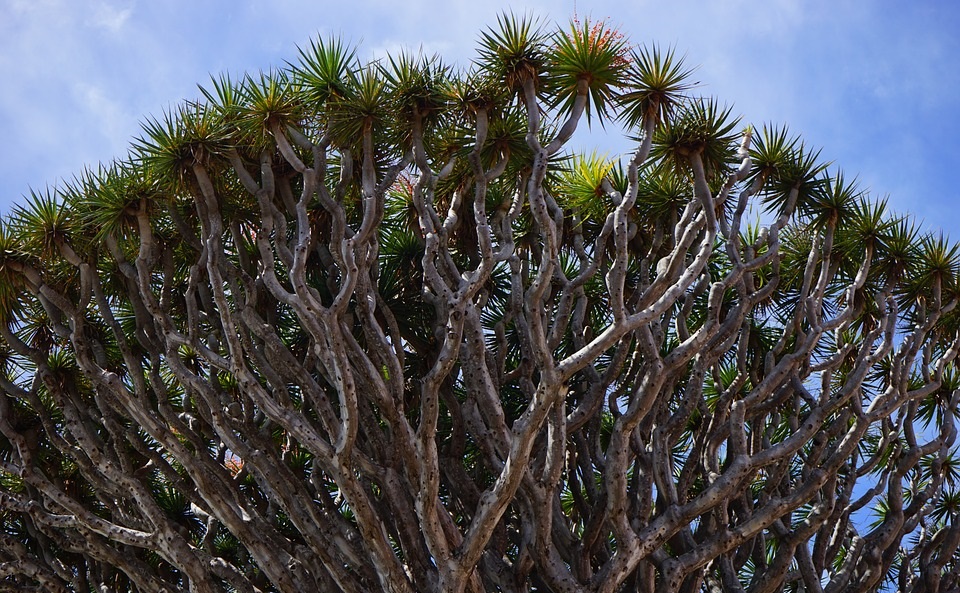}
\caption{Dracaena draco tree.}
\label{fig:draco}
\end{figure}

The goal of this paper is to provide an answer to this question when we focus on the vertical and horizontal displacements of a loaded binary tree. In particular, a binary tree perfectly characterizes a Dragon tree (Dracaena draco, see Fig. \ref{fig:draco}), which is a natural fractal-shaped structure and serves as a model. However, since the problem has been formulated abstractly, our results are general in nature and can be applied to all the structures that, whether natural or artificial (for example, pillars), meet the conditions of our model (Sec. \ref{sec:model}). There is a wide variety of methods to generate fractals with tree-like structure \cite{2} based on mathematical algorithms. Unfortunately, if the underlying physical mechanism responsible of such structure is unknown, then we cannot understand why fractals exist in Nature, and that is why we are so interested in obtaining a fractal by using natural laws and not algorithms.
\\

In nature, not only the laws of physics are present, but also the laws of biology. In a world where resources are scarce and energy is limited, the evolutionary process forces the geometry of structures to optimize resources. This is the underlying principle in Murray's Law \cite{3,4} that minimizes energy costs in blood transport in circulatory system and extends to respiratory system and to al vascular system in plants (Xylem). The minimum cost is achieved with the self-similar structure of the mentioned systems. Fractal structures designed to optimize performance are well known in nature and industry \cite{5,6} and there are more and more cases where fractals are responsible for optimization \cite{7,8,9,10}.

On the particular problem that concerns us, tree crowns are also subject to optimization processes. This can be seen in \cite{11}, when the maximum area for photosynthesis (minimun self-shading) with the minimum structural cost is required. Mandelbrot and Frame \cite{2,12} were the first to study crowns of binary self-similar trees into a plane, by using Iterated Function Systems (IFS) and geometric algorithms. However, in this paper, we will mathematically prove that the fractal structure of a tree crown emerges directly from the study of branching stresses, and for this we will simply use the continuum elasticity theory, not IFS; that is, the fractal structure of the tree crown emerges from natural principles and not from mathematical algorithms.

In addition, and this is an unexpected result, by using natural laws, we will find that fractals do not come alone, they can come in pairs. Generalizing this result, the question of which fractals and under which circumstances will appear in pairs arises. We are so used to generate natural-like fractals by using algorithms that we only pay attention to the fractals we are creating and not the others that are generated simultaneously. When we look at the crown of a tree, apparently we only see one fractal in its silhouette, but in reality there are two: the first one is associated to the vertical displacements of the branches and the second one to the horizontal displacements. In the first case, there is a Takagi curve, and in the second case, a linear combination of inverses of $\beta$-Cantor functions.

In order to achieve this result, we study (Sec. \ref{sec:model}) the well-known binary tree \cite{13} (see Fig. \ref{fig:struc}) that mimics the dragon tree. We will logically assume that the cross-section of the branches are smaller and smaller as we ascend in the successive levels of the tree as the stresses decrease. By applying the Principle of Virtual Work, we will calculate the vertical (Sec. \ref{sec:vertical}) and horizontal (Sec. \ref{sec:horizontal}) displacements in the $i$-th level of the tree. When the limit to infinity is taken, the Takagi curve (Sec. \ref{sec:vertical}) and a linear combination of inverses of $\beta$-Cantor functions (Sec. \ref{sec:horizontal}) give the vertical and horizontal displacements respectively. These two fractals are inescapably associated via the structure and their fractal dimension is linked as we prove in Sec. \ref{sec:dimensions}. This naturally brings to us the question of whether fractals will emerge in couples, which we will discuss in the conclusions (Sec. \ref{sec:conclusions}).

%% file: modelo.tex
\section{\bf{The MODEL: A STRUCTURE of $P$ LEVELS}}\label{sec:model}

\begin{figure}[h]
\centering
\includegraphics[width=9cm]{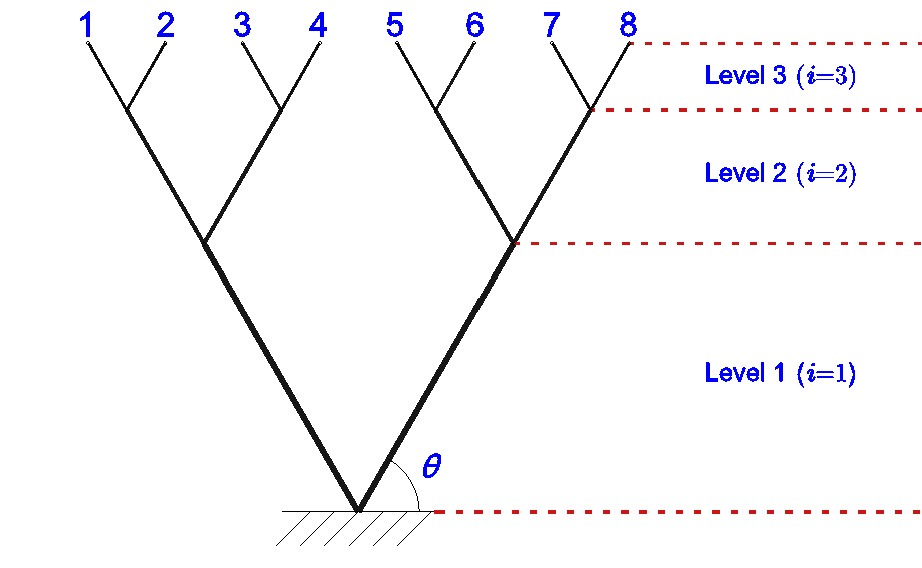}
\caption{Example of a binary tree structure with three levels, corresponding to $P=3$.}
\label{fig:struc}
\end{figure}

Let us consider a binary tree structure, with $P$ levels (see Fig. \ref{fig:struc}), such that the bifurcation points of the structure are rigid. The bars make an angle $\theta$ with respect to the horizontal axis and have the same Young’s modulus $E$ and shear modulus $G$ regardless of the level at which they are located. On the other hand, their length $L_i$ is determined by the level $i$ at which they are located and it is given by $L_i=L\cdot2^{1-i}, i=1,..,P$, being $L$ the length of the bars of the first level $i=1$.

The upper ends of the bars will be called nodes. Note that, except from the $P$-th level, the nodes coincide with the bifurcation points of the structure, where the lower bar forks generating two bars. The nodes of the $P$-th level are evenly distributed due to the ratio progression $1/2$ followed by the lengths of the bars from one level to another. These nodes of the $P$-th level (the higher one) will be called end nodes. Both nodes and bars will be listed from left to right, that is, $1,2,3,\ldots,2^i$ for a level $i$ (see Fig. \ref{fig:struc}).

The end nodes receive a total vertical load $F$ that is evenly distributed among them. We can interpret this load as the weight of the leaves on a tree whose structure is the given one and which corresponds to a Dracaena draco tree. Without loss of generality, we assume that the vertical load $F$ has a value of $1$. Therefore, each end node receive a downward vertical load of value $1/2^P$, corresponding to the weight of the leaves. 

The displacements of the nodes of the structure are given by the Principle of Virtual Work (PVW), which relates the loads and stresses in a virtual load system to the displacements and deformations in a real load system. The result of applying the PVW to our structure yields the following equation:
\begin{equation}\label{eq:1}
\centering
\sum_{i=1}^P \sum_{n=1}^{2^i} F_{i,n}^V\, \delta_{i,n}^R = \sum_{i=1}^P \sum_{j=1}^{2^i} \int_0^{L_i} \frac{M_{i,j}^V(x)\, M_{i,j}^R(x)}{E\, I_i } dx 
+\sum_{i=1}^P \sum_{j=1}^{2^i} \int_0^{L_i} \frac{N_{i,j}^V(x)\, N_{i,j}^R(x)}{E\, A_i } dx 
+\sum_{i=1}^P \sum_{j=1}^{2^i} \int_0^{L_i} \frac{T_{i,j}^V(x)\, T_{i,j}^R(x)}{G\, A_i^*} dx
\end{equation}

Where, for the $n$-th node, $n=1,\ldots,2^i$, of a level $i$:

\begin{itemize}
\item $F_{i,n}^V$: External force applied on the node on the virtual load system.
\item $\delta_{i,n}^R$: Node displacement on the real load system.
\end{itemize}

Where, for the $j$-th bar, $j=1,\ldots,2^i$, of a level $i$:

\begin{itemize}
\item $M_{i,j}^V(x)$: Bending moment on the virtual load system.
\item $M_{i,j}^R(x)$: Bending moment on the real load system.
\item $N_{i,j}^V(x)$: Axial stress on the virtual load system.
\item $N_{i,j}^R(x)$: Axial stress on the real load system.
\item $T_{i,j}^V(x)$: Shear stress on the virtual load system.
\item $T_{i,j}^R(x)$: Shear stress on the real load system.
\item $E$: Young's modulus of the bars.
\item $G$: Shear modulus of the bars.
\item $L_i$: Bar length.
\item $I$: Inertia of the bars on the first level $(i=1)$.
\item $a$: Ratio of inertia reduction per level, $a>1$.
\item $I_i$: Moment of inertia of the bar, $I_i=I\, a^{1-i}$.
\item $A$: Cross-sectional area of the bars on the first level $(i=1)$.
\item $u$: Ratio of cross-sectional area reduction per level, $u>1$.
\item $A_i$: Cross-sectional area of the bar, $A_i=A\, u^{1-i}$.
\item $A^*$: Cross-sectional shear area of the bars on the first level $(i=1)$.
\item $v$: Ratio of cross-sectional shear area reduction per level, $v>1$.
\item $A_i^*$: Cross-sectional shear area of the bar, $A_i^*=A^*\, v^{1-i}$.
\end{itemize}

As the axial and shear stresses are constant over the entire length of each bar, the Eq. (\ref{eq:1}) is rewritten as:
\begin{equation}\label{eq:2}
\sum_{i=1}^P \sum_{n=1}^{2^i} F_{i,n}^V\, \delta_{i,n}^R = \sum_{i=1}^P \sum_{j=1}^{2^i} \int_0^{L_i} \frac{M_{i,j}^V(x)\, M_{i,j}^R(x)}{E\, I_i } dx
+\sum_{i=1}^P \sum_{j=1}^{2^i} \frac{N_{i,j}^V\, N_{i,j}^R\, L\, 2^{1-i}}{E\, A_i } 
+\sum_{i=1}^P \sum_{j=1}^{2^i}\frac{T_{i,j}^V\, T_{i,j}^R\, L\, 2^{1-i}}{G\, A_i^* }
\end{equation}
where
\begin{equation}\label{eq:3}
\begin{gathered}
\sum_{i=1}^P \sum_{j=1}^{2^i} \int_0^{L_i} \frac{M_{i,j}^V(x)\, M_{i,j}^R(x)}{E\, I_i } dx 
\end{gathered}
\end{equation}

\begin{equation}\label{eq:4}
\begin{gathered}
\sum_{i=1}^P \sum_{j=1}^{2^i} \frac{N_{i,j}^V\, N_{i,j}^R\, L\, 2^{1-i}}{E\, A_i }
\end{gathered}
\end{equation}

\begin{equation}\label{eq:5}
\begin{gathered}
\sum_{i=1}^P \sum_{j=1}^{2^i}\frac{T_{i,j}^V\, T_{i,j}^R\, L\, 2^{1-i}}{G\, A_i^* }
\end{gathered}
\end{equation}
are the bending moments term (\ref{eq:3}), the axial stresses term (\ref{eq:4}) and the shear stresses term (\ref{eq:5}), respectively.

\begin{figure}[H]
\centering
\includegraphics[width=8cm]{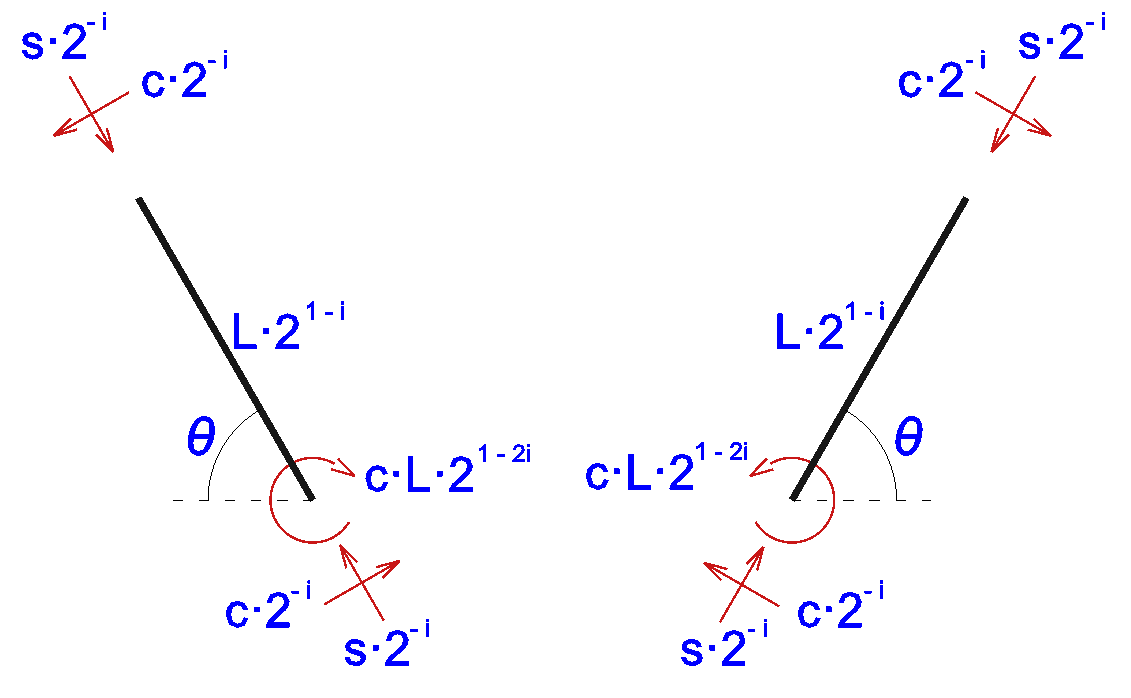}
\caption{Real load diagram of a bar of the $i$-th level (both possible positions). This bar has a length $L\, 2^{1-i}$ and makes an angle $\theta$ with respect to the horizontal. It is subjected to a bending moment $c\, L\, 2^{1-2i}$ in its bottom end and an axial stress $s\, 2^{-i}$ and shear stress $c\, 2^{-i}$ over its entire length.}
\label{fig:real}
\end{figure}

As a preliminary step to calculate the terms (\ref{eq:3}), (\ref{eq:4}) and (\ref{eq:5}), we first calculate the moments and stresses $M_{i,j}^R$, $N_{i,j}^R$, and $T_{i,j}^R$ produced by the real load system. Let us consider the $i$-th level of the structure, listing from the bottom up. The real load diagram of a bar of the $i$-th level is shown in Fig. \ref{fig:real}, considering that this bar can be positioned bottom-up and right-left or left-right:
\\
\noindent
where $c = \cos(\theta)$ and $s = \sin(\theta)$. Note that one bar of the $i$-th level connects to $2^{P-i}$ bars of the $P$-th level, and the loads of value $1/2^P$ are applied on the end nodes of these bars of the $P$-th level. Because of this, the real bending moments and axial and shear stresses are:
\begin{equation}\label{eq:6}
M_{i,j}^R(x) = c\, 2^{-i}(L\, 2^{1-i}-x), \, 0\leq x \leq L\, 2^{1-i}
\end{equation}
\begin{equation}\label{eq:7}
N_{i,j}^R = s\, 2^{-i}
\end{equation}
\begin{equation}\label{eq:8}
T_{i,j}^R = c\, 2^{-i}
\end{equation}

Note that on the right side of the PVW, the virtual and real stresses always appear multiplied together. Therefore, a sign criterion is not established for stresses in general, but the product will be considered as negative if the stresses have the opposite direction and positive if they have the same direction.

%% file: vertical.tex
\section{\bf{VERTICAL DISPLACEMENTS and the TAKAGI CURVE}}\label{sec:vertical}
\subsection{\bf{Vertical Displacements in the Finite Structure}}\label{sec:vfinite}

As a previous step to calculate the vertical displacements in the structure with infinite levels, is necessary to calculate the vertical displacements of the end nodes in a finite structure. For this, we will calculate the displacement $V\delta_{P,w}^R$ of an arbitrary end node $w$ in a structure with $P$ levels by using the Eq. (\ref{eq:2}) on a virtual load system consisting of a downward vertical load of value $1$ on the end node $w$. Recall that $V\delta_{P,w}^R$ are the vertical displacements per unit load, as the real load is $F=1$.

The load is transmitted from a bar to the one connected on the lower level and so on until it reaches the base. This will generate at the base a vertical stress, of value $1$, plus a bending moment. Note that there will only be one loaded bar per level. In addition, as there is no external load on any node of the structure except for the end node $w$, the left side of Eq. (\ref{eq:2}) is rewritten as:
\begin{equation}\label{eq:9}
\sum_{i=1}^P \sum_{n=1}^{2^i} F_{i,n}^V\, \delta_{i,n}^R = 1\, V\delta_{P,w}^R = V\delta_{P,w}^R
\end{equation}

\begin{figure}[h]
\centering
\includegraphics[width=9cm]{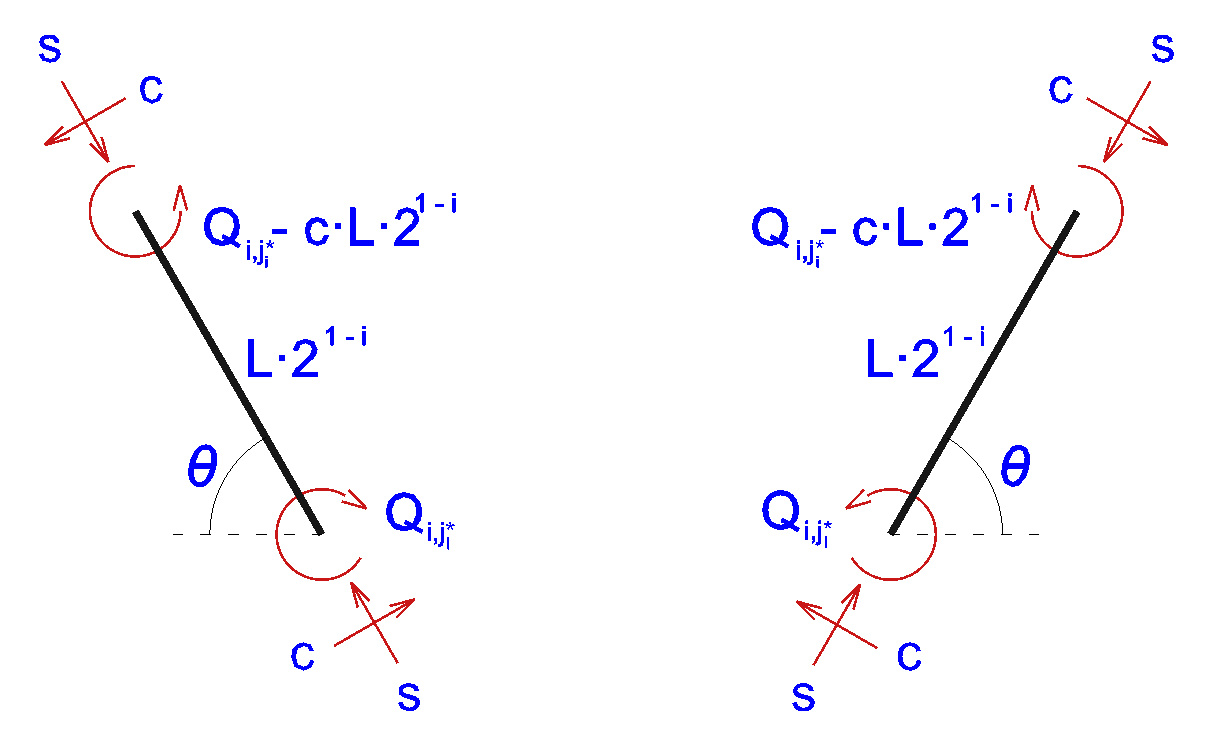}
\caption{Virtual load diagram of the loaded bar $j_i^*$ of the $i$-th level for vertical displacements (both possible positions). This bar has a length $L\, 2^{1-i}$ and makes an angle $\theta$ with respect to the horizontal. It is subjected to a bending moment $Q_{i,j_i^*}$ in its bottom end and an axial stress $s$ and shear stress $c$ over its entire length.}
\label{fig:vvirtual}
\end{figure}

The calculation of the terms (\ref{eq:3}), (\ref{eq:4}) and (\ref{eq:5}), on the right side of Eq. (\ref{eq:2}), requires calculating the stresses $M_{i,j}^V$, $N_{i,j}^V$, and $T_{i,j}^V$ produced by the virtual load system. As already indicated, there is only one loaded bar per level. By naming $j_i^*$ that bar for the $i$-th level, its loading diagram is shown in Fig. \ref{fig:vvirtual}, considering that this bar can be positioned bottom-up and right-left or left-right.
\\

Therefore, the stresses on the bar $j_i^*$ are:
\begin{equation}\label{eq:10}
M_{i,j}^V(x) = Q_{i,j_i^*}-c\, x, 0\leq x \leq L\, 2^{1-i}
\end{equation}
\begin{equation}\label{eq:11}
N_{i,j}^V = s
\end{equation}
\begin{equation}\label{eq:12}
T_{i,j}^V = c
\end{equation}
where $Q_{i,j_i^*}$ is the moment generated at the bottom end of the loaded bar $j_i^*$. To calculate this moment, we need to calculate the horizontal distance $DH_{i,j_i^*}(w)$ between the bottom end of the loaded bar $j_i^*$ and the end node $w$ where the load is applied. Let us assume without loss of generality that the end nodes are uniformly distributed in the interval $[0,1]$, such that the end node $1$ is located at $z=\left( \frac{1}{2} \right)^{P+1}$ in that interval and the end node $2^P$ at $z=1-\left(\frac{1}{2} \right)^{P+1}$. As the nodes are equispaced, it results that the end node $w$ is located at:
\begin{equation}\label{eq:13}
z(w)=\frac{w-1}{2^P-1}\left(1-\frac{1}{2^P}\right)+\frac{1}{2^{P+1}}
\end{equation}
As a consequence of setting the position of the end nodes, the length $L$ and the tilt angle $\theta$ are related.

According to Fig. \ref{fig:struc}, and taking into account the geometrical progression of ratio $1/2$ followed by the lengths of the bars, the horizontal distance $DH_{i,j_i^*}(w)$ between the bottom end of the loaded bar $j_i^*$ and the end node $w$ where the load is applied is:
\begin{equation}\label{eq:14}
DH_{i,j_i^*}(w)=4\, c\, L\left(\frac{1}{2^i}-\frac{\sigma \left(2^{i-1}\, z(w)\right)}{2^{i-1}}\right)
\end{equation}

Where $\sigma(x)=\min_{k\in \mathbb{N}} |x-k|$.

Due to the geometry of the structure, the virtual moment generated in the base $Q_{i,j_i^*}$ will always coincide in direction with the real moment generated in the base. Therefore, the moment $Q_{i,j_i^*}$ caused by the unit load of the virtual system is:
\begin{equation}\label{eq:15}
Q_{i,j_i^*}=1\, DH_{i,j_i^*}(w)=4\, c\, L\left(\frac{1}{2^i}-\frac{\sigma \left(2^{i-1}\, z(w)\right)}{2^{i-1}}\right)=
\frac{c\, L}{2^{i-2}}-\frac{c\, L\, \sigma \left(2^{i-1}\, z(w)\right)}{2^{i-3}}
\end{equation}

\subsubsection{\bf{Bending Moments}}\label{sec:vmoment}

As the bar $j_i^*$ it is the only one loaded on the $i$-th level and $E$ and $I_i=I\,a^{1-i}$ are constants for each bar, the bending moments term (\ref{eq:3}), according to (\ref{eq:6}), (\ref{eq:10}) and (\ref{eq:15}), is rewritten as:
\begin{equation}
\sum_{i=1}^{P} \frac{1}{E I_{i}} \int_{0}^{L_{i}} M_{i, j_{i}^{*}}^{V}(x) \, M_{i, j_{i}^{*}}^{R}(x) \, d x=\frac{20 \, c^{2} \, L^{3}}{3 \, E \, I \, a} \sum_{i=1}^{P}\left(\frac{a}{16}\right)^{i}-\frac{c^{2} \, L^{3}}{E \, I} \sum_{i=1}^{P}\left(\frac{a}{16}\right)^{i-1} \sigma\left(2^{i-1} z(w)\right)=\nonumber
\end{equation}
\begin{subequations}
\begin{empheq}[left={=\empheqlbrace}]{align} 
\begin{aligned}[t] \label{eq:18a}
           & \frac{20 \, c^{2} \, L^{3}}{3 \, E \, I \, a}\frac{1-\left(\frac{a}{16}\right)^{P}}{\frac{16}{a}-1}-\frac{c^{2} \, L^{3}}{E \, I} \phi\left(z(w) ; \frac{a}{16}, P-1\right) & \text {if} \quad a \neq 16
           \end{aligned}\\
\begin{aligned}[t] \label{eq:18b}
	    & \frac{5 \, c^{2} \, L^{3} \, P}{12 \, E \, I}-\frac{c^{2} \, L^{3}}{E \, I} \, \phi(z(w) ; 1, P-1) & \text {if} \quad a=16
	   \end{aligned}
\end{empheq}
\end{subequations}
where
\begin{equation} \label{eq:19}
\phi\left(z(w) ; \frac{a}{16}, P-1\right)=\sum_{i=1}^{P}\left(\frac{a}{16}\right)^{i-1} \sigma\left(2^{i-1} z(w)\right)
\end{equation}
\\
\subsubsection{\bf{Axial Stresses}}\label{sec:vaxial}

To calculate the axial stresses we proceed in a similar way to that used for the calculation of the bending moments (see Sec. \ref{sec:vmoment}). As the bar $j_i^*$ is the only one loaded on the $i$-th level and $A_i = A\, u^{1-i}$, the axial stresses term (\ref{eq:4}), according to (\ref{eq:7}) and (\ref{eq:11}), is rewritten as:

\begin{equation}
\sum_{i=1}^{P} \frac{N_{i, j_{i}^{*}}^{V} \, N_{i, j_{i}^{*}}^{R} \, L \, 2^{1-i}}{E \, A_{i}} 
=\sum_{i=1}^{P} \frac{s^{2} \, 2^{-i} \, L \, 2^{1-i}}{E \, A \, u^{1-i}}=
\frac{2 \, s^{2} \, L}{E \, A \, u} \sum_{i=1}^{P}\left(\frac{u}{4}\right)^{i}= \nonumber
\end{equation}

\begin{subequations}
\begin{empheq}[left={=\empheqlbrace}]{align}
\begin{aligned}[t] \label{eq:21a}
           & \frac{2 \, s^{2} \, L}{E \, A \, u}\frac{1-\left(\frac{u}{4}\right)^{P}}{\frac{4}{u}-1}\quad \text {if} \quad u \neq 4
           \end{aligned}\\
\begin{aligned}[t] \label{eq:21b}
	    & \frac{s^{2} \, L \, P}{2 \, E \, A} \quad \text {if} \quad u=4
	   \end{aligned}
\end{empheq}
\end{subequations}
\\
\subsubsection{\bf{Shear Stresses}}\label{sec:vshear}

We proceed similarly to the previous cases (see Secs. \ref{sec:vmoment} and \ref{sec:vaxial}). As the bar $j_i^*$ is the only one loaded on the $i$-th level and $A_i^* = A^*\, v^{1-i}$, the shear stresses term (\ref{eq:5}), according to (\ref{eq:8}) and (\ref{eq:12}), is rewritten as:

\begin{equation}
\sum_{i=1}^{P} \frac{T_{i, j_{i}^{*}}^{V} \, T_{i, j_{i}^{*}}^{R} \, L \, 2^{1-i}}{G \, A_{i}^*}
=\sum_{i=1}^{P} \frac{c^{2} \, 2^{-i} \, L \, 2^{1-i}}{G \, A^* \, v^{1-i}}=
\frac{2 \, c^{2} \, L}{G \, A^* \, v} \sum_{i=1}^{P}\left(\frac{v}{4}\right)^{i}= \nonumber
\end{equation}

\begin{subequations}
\begin{empheq}[left={=\empheqlbrace}]{align}
\begin{aligned}[t] \label{eq:23a}
           & \frac{2 \, c^{2} \, L}{G \, A^* \, v}\frac{1-\left(\frac{v}{4}\right)^{P}}{\frac{4}{v}-1}\quad \text {if} \quad v \neq 4
           \end{aligned}\\
\begin{aligned}[t] \label{eq:23b}
	    & \frac{c^{2} \, L \, P}{2 \, G \, A^*} \quad \text {if} \quad v=4
	   \end{aligned}
\end{empheq}
\end{subequations}
\\
\subsubsection{\bf{Total Displacements in the Finite Structure}}\label{sec:vtotal}

The total vertical displacement $V\delta_{P,w}^R$ of the end node located at $z(w)$ is obtained by adding up the expressions (\ref{eq:18a}), (\ref{eq:18b}), (\ref{eq:21a}), (\ref{eq:21b}), (\ref{eq:23a}) and (\ref{eq:23b}), depending on the different values of the parameters $a$, $u$ and $v$.

\begin{widetext}
\begin{itemize}
\item For $a\neq16$, $u\neq4$, $v\neq4$ (terms (\ref{eq:18a}), (\ref{eq:21a}) and (\ref{eq:23a}):
\begin{equation} \label{eq:24}
V \delta_{P, w}^{R}=\frac{20 \, c^{2} \, L^{3}}{3 \, E \, I \, a} \frac{1-\left(\frac{a}{16}\right)^{P}}{\frac{16}{a}-1}-\frac{c^{2} \, L^{3}}{E \, I}\phi\left(z(w) ; \frac{a}{16}, P-1\right)+\frac{2 \, s^{2} \, L}{E \, A \, u} \frac{1-\left(\frac{u}{4}\right)^{P}}{\frac{4}{u}-1}+\frac{2 \, c^{2} \, L}{G \, A^{*} \, v} \frac{1-\left(\frac{v}{4}\right)^{P}}{\frac{4}{v}-1}
\end{equation}

\item For $a\neq16$, $u\neq4$, $v=4$ (terms (\ref{eq:18a}), (\ref{eq:21a}) and (\ref{eq:23b}):
\begin{equation} \label{eq:25}
V \delta_{P, w}^{R}=\frac{20 \, c^{2} \, L^{3}}{3 \, E \, I \, a} \frac{1-\left(\frac{a}{16}\right)^{P}}{\frac{16}{a}-1}-\frac{c^{2} \, L^{3}}{E \, I}\phi\left(z(w) ; \frac{a}{16}, P-1\right)+\frac{2 \, s^{2} \, L}{E \, A \, u} \frac{1-\left(\frac{u}{4}\right)^{P}}{\frac{4}{u}-1}+\frac{c^{2} \, L \, P}{2 \, G \, A^{*}}
\end{equation}

\item For $a\neq16$, $u=4$, $v\neq4$ (terms (\ref{eq:18a}), (\ref{eq:21b}) and (\ref{eq:23a}):
\begin{equation} \label{eq:26}
V \delta_{P, w}^{R}=\frac{20 \, c^{2} \, L^{3}}{3 \, E \, I \, a} \frac{1-\left(\frac{a}{16}\right)^{P}}{\frac{16}{a}-1}-\frac{c^{2} \, L^{3}}{E \, I}\phi\left(z(w) ; \frac{a}{16}, P-1\right)+\frac{s^{2} \, L \, P}{2 \, E \, A}+\frac{2 \, c^{2} \, L}{G \, A^{*} \, v} \frac{1-\left(\frac{v}{4}\right)^{P}}{\frac{4}{v}-1}
\end{equation}

\item For $a\neq16$, $u=4$, $v=4$ (terms (\ref{eq:18a}), (\ref{eq:21b}) and (\ref{eq:23b}):
\begin{equation} \label{eq:27}
V \delta_{P, w}^{R}=\frac{20 \, c^{2} \, L^{3}}{3 \, E \, I \, a} \frac{1-\left(\frac{a}{16}\right)^{P}}{\frac{16}{a}-1}-\frac{c^{2} \, L^{3}}{E \, I}\phi\left(z(w) ; \frac{a}{16}, P-1\right)+\frac{s^{2} \, L \, P}{2 \, E \, A}+\frac{c^{2} \, L \, P}{2 \, G \, A^{*}}
\end{equation}

\item For $a=16$, $u\neq4$, $v\neq4$ (terms (\ref{eq:18b}), (\ref{eq:21a}) and (\ref{eq:23a}):
\begin{equation} \label{eq:28}
V \delta_{P, w}^{R}=\frac{5 \, c^{2} \, L^{3} \, P}{12 \, E \, I \, a}-\frac{c^{2} \, L^{3}}{E \, I}\phi\left(z(w) ; 1, P-1\right)+\frac{2 \, s^{2} \, L}{E \, A \, u} \frac{1-\left(\frac{u}{4}\right)^{P}}{\frac{4}{u}-1}+\frac{2 \, c^{2} \, L}{G \, A^{*} \, v} \frac{1-\left(\frac{v}{4}\right)^{P}}{\frac{4}{v}-1}
\end{equation}

\item For $a=16$, $u\neq4$, $v=4$ (terms (\ref{eq:18b}), (\ref{eq:21a}) and (\ref{eq:23b}):
\begin{equation} \label{eq:29}
V \delta_{P, w}^{R}=\frac{5 \, c^{2} \, L^{3} \, P}{12 \, E \, I \, a}-\frac{c^{2} \, L^{3}}{E \, I}\phi\left(z(w) ; 1, P-1\right)+\frac{2 \, s^{2} \, L}{E \, A \, u} \frac{1-\left(\frac{u}{4}\right)^{P}}{\frac{4}{u}-1}+\frac{c^{2} \, L \, P}{2 \, G \, A^{*}}
\end{equation}

\item For $a=16$, $u=4$, $v\neq4$ (terms (\ref{eq:18b}), (\ref{eq:21b}) and (\ref{eq:23a}):
\begin{equation} \label{eq:30}
V \delta_{P, w}^{R}=\frac{5 \, c^{2} \, L^{3} \, P}{12 \, E \, I \, a}-\frac{c^{2} \, L^{3}}{E \, I}\phi\left(z(w) ; 1, P-1\right)+\frac{s^{2} \, L \, P}{2 \, E \, A}+\frac{2 \, c^{2} \, L}{G \, A^{*} \, v} \frac{1-\left(\frac{v}{4}\right)^{P}}{\frac{4}{v}-1}
\end{equation}

\item For $a=16$, $u=4$, $v=4$ (terms (\ref{eq:18b}), (\ref{eq:21b}) and (\ref{eq:23b}):
\begin{equation} \label{eq:31}
V \delta_{P, w}^{R}=\frac{5 \, c^{2} \, L^{3} \, P}{12 \, E \, I \, a}-\frac{c^{2} \, L^{3}}{E \, I}\phi\left(z(w) ; 1, P-1\right)+\frac{s^{2} \, L \, P}{2 \, E \, A}+\frac{c^{2} \, L \, P}{2 \, G \, A^*}
\end{equation}
\end{itemize}
\end{widetext}

\subsection{\bf{Vertical Displacements in the Structure with Infinite Levels: The Takagi Curve}}

To study the structure with infinite levels, the limit $P\rightarrow\infty$ is taken. As a result, the expression (\ref{eq:19}) becomes the Takagi curve:

\begin{equation}
\lim _{P \rightarrow \infty} \phi\left(z(w) ; \frac{a}{16}, P-1\right)=\Psi_{\frac{a}{16}}(z(w)) \quad \text { with }|a|<16
\end{equation}

\begin{figure}[h]
\centering
\includegraphics[width=10cm]{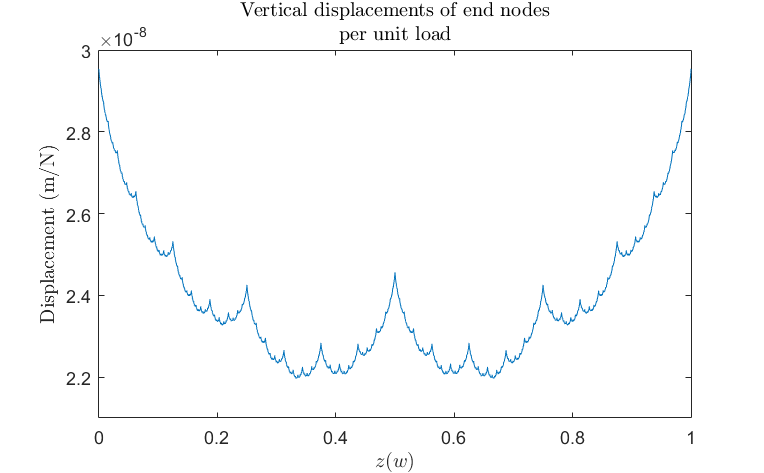}
\caption{Takagi curve for the vertical displacements of the end nodes per unit load versus the end nodes positions $z(w)$ for a structure with infinite levels. A positive value indicates downward displacement. The vertical displacements are plotted for $\theta=60^{\circ}$, $E=10^{10}$ N/m\textsuperscript{2}, $G=5\cdot10^8$ N/m\textsuperscript{2}, $L=0.5$ m, $I=3.1416\cdot10^{-4}$ m\textsuperscript{4}, $A=3.1416\cdot10^{-2}$ m\textsuperscript{2}, $A^*=2.8274\cdot10^{-2}$ m\textsuperscript{2}, $a=9$, $u=3$ and $v=3$.}
\label{fig:vdesp}
\end{figure}

Where $\Psi_\frac{a}{16}\left(z\left(w\right)\right)$ belongs to the exponential Takagi class \cite{20}. Due to mechanical conditions, $a>1$, so the Takagi curve would be obtained for values  $1< a<16$ (see Fig. \ref{fig:vdesp}). To study the vertical displacements at the end nodes when $P\rightarrow\infty$, we must take limit $P\rightarrow\infty$ on Eqs. (\ref{eq:24})-(\ref{eq:31}), resulting in:

\begin{itemize}
\item For $1< a<16$, $1<u<4$ and $1<v<4$:
\begin{equation} \label{eq:33}
V \delta_{P\rightarrow\infty, w}^{R}=\frac{20 \, c^{2} \, L^{3}}{3 \, E \, I(16-a)} 
-\frac{c^{2} \, L^{3}}{E \, I}\Psi_{\frac{a}{16}}(z(w))
+\frac{2 \, s^{2} \, L}{E \, A(4-u)}+\frac{2 \, c^{2} \, L}{G \, A^{*}(4-v)}
\end{equation}

\item For $1<a<16$, $u\geq4$ and/or $v\geq4$, $V \delta_{P\rightarrow\infty, w}^{R}\rightarrow\infty$ for all $w$ according to a geometric series of ratio greater than or equal to one (since the Takagi curve is bounded). Consequently, the structure collapses.
\end{itemize}

%% file: horizontal.tex
\section{\bf{HORIZONTAL DISPLACEMENTS and the $\beta$-CANTOR FUNCTION}}\label{sec:horizontal}
\subsection{\bf{Horizontal Displacements in the Finite Structure}}\label{sec:hfinite}

To calculate the horizontal displacements $H\delta_{P,w}^R$ of the end nodes using Eq. (\ref{eq:2}), we proceed in a similar way to that used in Sec. \ref{sec:vertical}. Note that the structure and the real load system are symmetrical with respect to an axis perpendicular to the ground passing through the soil-structure junction point. Thus, due to ease of calculation, only the left half of the strucure will be considered. We treat a virtual load system consisting of a leftward horizontal load of value $1$ on the end node $w$ whose displacement is to be calculated. This load is transmitted from a bar to the one connected on the lower level and so on until it reaches the base. This will generate a horizontal  stress of value $1$ plus a moment at the base. Recall that there will only be one bar loaded per level. Recall that $H\delta_{P,w}^R$ are the horizontal displacements per unit load, as the real load is $F=1$.

\begin{figure}[h]
\centering
\includegraphics[width=9cm]{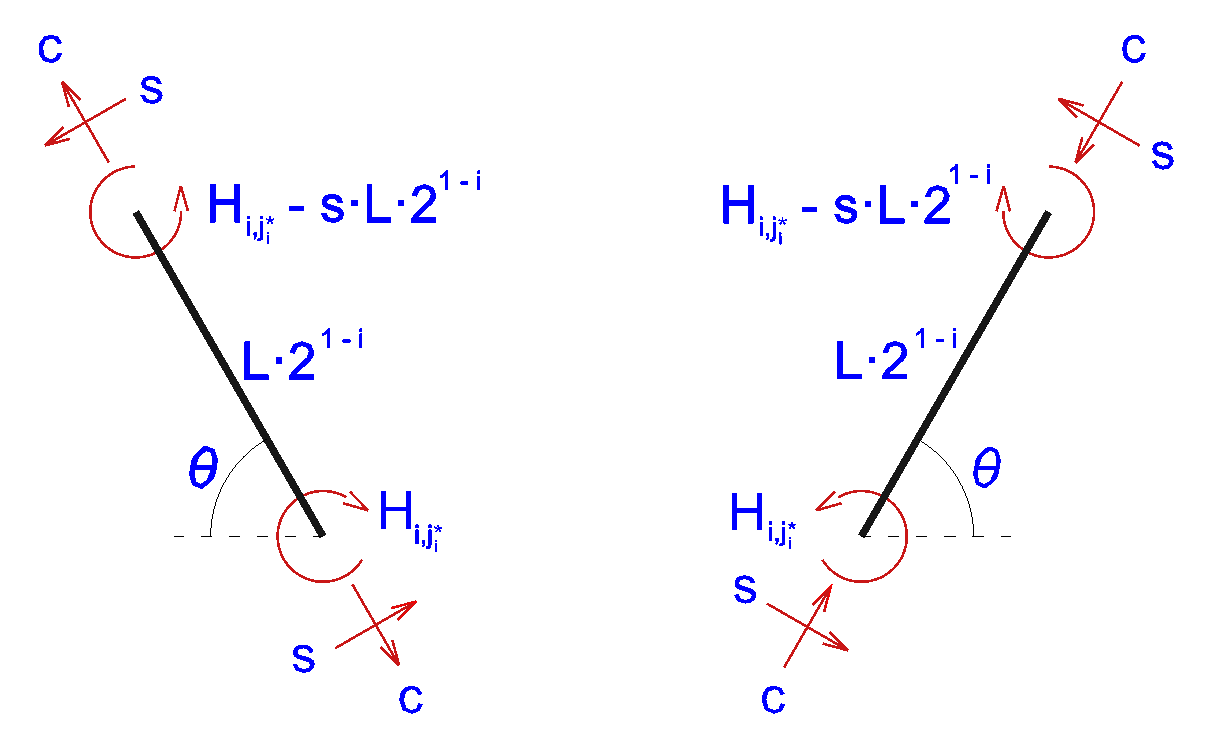}
\caption{Virtual load diagram of the loaded bar $j_i^*$ of the $i$-th level for horizontal displacements (both possible positions). This bar has a length $L\, 2^{1-i}$ and makes an angle $\theta$ with respect to the horizontal. It is subjected to a bending moment $H_{i,j_i^*}$ in its bottom end and an axial stress $c$ and shear stress $s$ over its entire length.}
\label{fig:hvirtual}
\end{figure}

To calculate the terms (\ref{eq:3}), (\ref{eq:4}) and (\ref{eq:5}) on the right side of Eq. (\ref{eq:2}) we proceed in a similar way to Sec. \ref{sec:vertical}. To this end, we must calculate the stresses $M_{i,j}^V$, $P_{i,j}^V$, and $T_{i,j}^V$ as shown in the virtual load diagram of the loaded bar $j_i^*$ on Fig. \ref{fig:hvirtual}.
\\

Without loss of generality, let us suppose that the end nodes are uniformly distributed in the interval $[0,1],$ with node $1$ located at $z^\ast=0$ and the node $2^P$ at $z^\ast=1$. Thus, as the nodes are uniformly distributed, the end node $w$ is located at:
\begin{equation}
z^{*}(w)=\frac{w-1}{2^{P}-1} \nonumber
\end{equation}

Due to ease of calculation, these positions $z^*(w)$ are different from $z(w)$ of Sec. \ref{sec:vfinite}, but coincide when the limit $P\rightarrow\infty$ is taken. According to Fig. \ref{fig:hvirtual}, the stresses of the loaded bar $j_i^*$ are:
\begin{equation} \label{eq:34}
M_{i, j t}^{V}(x)=\left(H_{i, j_{i}^{*}}-s \, x\right)\left(1-2 \rho_{i}\left(\sigma\left(z^{*}(w)\right)\right)\right), \,
0 \leq x \leq L \, 2^{1-i}
\end{equation}
\begin{equation} \label{eq:35}
N_{i, j_{i}^{*}}^{V}=c\left(2 \rho_{i}\left(\sigma\left(z^{*}(w)\right)\right)-1\right)
\end{equation}
\begin{equation} \label{eq:36}
T_{i, j i}^{V}=s\left(1-2 \rho_{i}\left(\sigma\left(z^{*}(w)\right)\right)\right)
\end{equation}
where $\rho_i(x)\in\left\{0,1\right\}$ is the $i$-th coefficient of the dyadic expansion of $x$, that is:
\begin{equation} \label{eq:37}
x=\sum_{k=1}^{\infty} \frac{\rho_{k}(x)}{2^{k}}
\end{equation}

As $\rho_i(x)$ only takes the values $0$ or $1$, the expression $\left(2\rho_i(x)-1\right)$ takes values $-1$ and $+1$ respectively. Therefore, this expression indicates the sign of the stresses on the bar: $-1$ if the bar is positioned bottom-up and right-left and $+1$ if the bar is positioned bottom-up and left-right.

On the other hand, $H_{i,j_i^\ast}$ is the virtual moment generated in the bottom end of the bar. To calculate this moment, we need to know the vertical distance $DV_{i,j_i^\ast}(w)$ from the bottom end of the loaded bar $j_i^\ast$ to the loaded end node $w$. This distance $DV_{i,j_i^\ast}(w)$ is:

\begin{equation} \label{eq:38}
D V_{i, j_{i}^{*}}(w)=s \, L\left(\left(\frac{1}{2}\right)^{i-2}-\left(\frac{1}{2}\right)^{P-1}\right)
\end{equation}

Therefore, the moment $H_{i,j_i^\ast}$ obtained from Eq. (\ref{eq:38}) is:

\begin{equation} \label{eq:39}
H_{i, j_{i}^{*}}=1 \, D V_{i, j_{i}^{*}}(w)=s \, L\left(\left(\frac{1}{2}\right)^{i-2}-\left(\frac{1}{2}\right)^{P-1}\right)
\end{equation}

\subsubsection{\bf{Bending Moments}}\label{sec:hmoment}

We proceed in a similar way as in Sec. \ref{sec:vmoment}. As the bar $j_i^*$ is the only one loaded on the $i$-th level and $E$ and $I_i=I\, a^{1-i}$ are constants for each bar, the bending moments term (\ref{eq:3}), according to (\ref{eq:6}), (\ref{eq:34}) and (\ref{eq:39}), is rewritten as:
\begin{equation} \label{eq:42}
\begin{gathered}
\sum_{i=1}^{P} \frac{1}{E I_{i}} \int_{0}^{L_{i}} M_{i, j_{i}^{*}}^{V}(x) \, M_{i, j_{i}^{*}}^{R}(x) \, d x=\\
=\frac{4 \, c \, s \, L^{3}}{E \, I \, a}\left(\frac{5}{3}\sum_{i=1}^{P}\left(\frac{a}{16}\right)^{i}-\frac{10}{3} \sum_{i=1}^{P} \rho_{i}\left(\sigma\left(z^{*}(w)\right)\right)\left(\frac{a}{16}\right)^{i}-\frac{1}{2^{P}} \sum_{i=1}^{P}\left(\frac{a}{8}\right)^{i}+\frac{1}{2^{P-1}} \sum_{i=1}^{P} \rho_{i}\left(\sigma\left(z^{*}(w)\right)\right)\left(\frac{a}{8}\right)^{i}\right)=
\end{gathered} 
\end{equation}
\begin{subequations}
\begin{empheq}[left={=\empheqlbrace}]{align}
\begin{aligned}[t] \label{eq:43a}
	    &\frac{4 \, c \, s \, L^{3}}{E \, I}\left(\frac{5-5\left(\frac{a}{16}\right)^{P}}{48-3\,a}-\frac{10}{3 \, a} C_{P, \frac{a}{16}}\left(\sigma\left(z^{*}(w)\right)\right)+\right.\\
	    &\left.-\frac{\left(\frac{1}{2}\right)^P-\left(\frac{a}{16}\right)^P}{8-a}+\frac{1}{2^{P-1} \, a} C_{P,\frac{a}{8}}\left(\sigma\left(z^{*}(w)\right)\right)\right) & \text {if} \quad a\neq8, a\neq16
           \end{aligned}\\
\begin{aligned}[t] \label{eq:43b}
	    &\frac{4 \, c \, s \, L^{3}}{E \, I}\left(\frac{5}{3} \frac{1-1/16^{P}}{8}-\frac{5}{12} C_{P, \frac{1}{2}}\left(\sigma\left(z^{*}(w)\right)\right)-\frac{P}{2^{P+3}}-\frac{1}{2^{P+2}} C_{P,1}\left(\sigma\left(z^{*}(w)\right)\right)\right) & \text {if} \quad a=8
           \end{aligned}\\
\begin{aligned}[t] \label{eq:43c}
	    &\frac{4 \, c \, s \, L^{3}}{E \, I}\left(\frac{5\, P}{48}-\frac{5}{24} C_{P, 1}\left(\sigma\left(z^{*}(w)\right)\right)+\frac{1/2^P-1}{8}+\frac{1}{2^{P+3}} C_{P,2}\left(\sigma\left(z^{*}(w)\right)\right)\right) & \text {if} \quad a=16
           \end{aligned}
\end{empheq}
\end{subequations}
Where
\begin{equation} \label{eq:44}
C_{P, t}(x)=\sum_{k=1}^{P} \rho_{k}(x) \, t^{k}
\end{equation}
\\
\subsubsection{\bf{Axial Stresses}}\label{sec:haxial}

We proceed in a similar way as in Sec. \ref{sec:vaxial}. The axial stresses term is given by (\ref{eq:4}). Therefore, substituting $A_i=A\, u^{1-i}$, (\ref{eq:7}) and (\ref{eq:35}) into (\ref{eq:4}), it results:

\begin{equation} \label{eq:45}
\sum_{i=1}^{P} \frac{N_{i, j_{i}^{*}}^{V} \, N_{i, j_{i}^{*}}^{R} \, L \, 2^{1-i}}{E \, A_{i}}=\frac{2 \, c \, s \, L}{E \, A \, u} \sum_{i=1}^{P}\left(2 \rho_{i}\left(\sigma\left(z^{*}(w)\right)\right)-1\right) \frac{u^{i}}{4^{i}}=
\end{equation}
\begin{subequations}
\begin{empheq}[left={=\empheqlbrace}]{align}
\begin{aligned}[t] \label{eq:46a}
	    &\frac{2 \, c \, s \, L}{E \, A \, u}\left(2 \, C_{P, \frac{u}{4}}\left(\sigma\left(z^{*}(w)\right)\right)-\frac{1-\left(\frac{u}{4}\right)^{P}}{\frac{4}{u}-1}\right) & \text {if} \quad u\neq4
           \end{aligned}\\
\begin{aligned}[t] \label{eq:46b}
	    &\frac{c \, s \, L}{2\, E \, A}\left(2 \, C_{P, 1}\left(\sigma\left(z^{*}(w)\right)\right)\right) & \text {if} \quad u=4
           \end{aligned}
\end{empheq}
\end{subequations}
\\
\subsubsection{\bf{Shear Stresses}}\label{sec:hshear}

We proceed in a similar way as in Sec. \ref{sec:vshear}. The shear stresses term is given by (\ref{eq:5}). Therefore, substituting $A^*_i=A^*\, v^{1-i}$, (\ref{eq:8}) and (\ref{eq:36}) into (\ref{eq:5}), it results:

\begin{equation} \label{eq:47}
\sum_{i=1}^{P} \frac{T_{i, j_{i}^{*}}^{V} \, T_{i, j_{i}^{*}}^{R} \, L \, 2^{1-i}}{G \, A_{i}^*}=\frac{2 \, c \, s \, L}{G \, A^* \, v} \sum_{i=1}^{P}\left(1-2 \rho_{i}\left(\sigma\left(z^{*}(w)\right)\right)\right) \frac{v^{i}}{4^{i}}=
\end{equation}
\begin{subequations}
\begin{empheq}[left={=\empheqlbrace}]{align}
\begin{aligned}[t] \label{eq:48a}
	    &\frac{2 \, c \, s \, L}{G \, A^{*} \, v}\left(\frac{1-\left(\frac{v}{4}\right)^{P}}{\frac{4}{v}-1}-2 \, C_{P, \frac{v}{4}}\left(\sigma\left(z^{*}(w)\right)\right)\right) & \text {if} \quad v\neq4
           \end{aligned}\\
\begin{aligned}[t] \label{eq:48b}
	    &\frac{c \, s \, L}{2\, G \, A^{*}}\left(P-2 \, C_{P, 1}\left(\sigma\left(z^{*}(w)\right)\right)\right) & \text {if} \quad v=4
           \end{aligned}
\end{empheq}
\end{subequations}
\\
\subsubsection{\bf{Total Displacements in the Finite Structure}}\label{sec:htotal}

The horizontal displacement $H\delta_{P,w}^R$ of the node $w$ located at $z^\ast(w)$ is obtained by adding up the expressions (\ref{eq:43a}), (\ref{eq:43b}), (\ref{eq:43c}), (\ref{eq:46a}), (\ref{eq:46b}), (\ref{eq:48a}) and (\ref{eq:48b}) depending on the parameters $a$, $u$ and $v$. The displacement is to the left if the node $w$ belongs to the set $\{1,\ldots,2^{P-1}\}$ and to the right if it belongs to the set $\{2^{P-1}+1,\ldots,2^P\}$.

\begin{itemize}
\item For $a\neq16$, $a\neq8$, $u\neq4$ and $v\neq4$ (terms (\ref{eq:43a}), (\ref{eq:46a}) and (\ref{eq:48a})):
\begin{equation} \label{eq:49}
\begin{gathered}
H \delta_{P, w}^{R}=2 \, c \, s \, L\left(\frac{10 \, L^{2}}{3 \, E \, I} \, \frac{1-\left(\frac{a}{16}\right)^{P}}{16-a}-\frac{20 \, L^{2}}{3 \, E \, I \, a} C_{P, \frac{a}{16}}\left(\sigma\left(z^{*}(w)\right)\right)-2 \, L^{2} \frac{\left(\frac{1}{2}\right)^P-\left(\frac{a}{16}\right)^{P}}{E \, I(8-a)}+\right.\\
+\frac{L^{2}}{2^{P-2} \, E \, I \, a} C_{P, \frac{a}{8}}\left(\sigma\left(z^{*}(w)\right)\right)+ \frac{2}{E \, A \, u} C_{P, \frac{u}{4}}\left(\sigma\left(z^{*}(w)\right)\right)-\frac{1}{E \, A \, u} \frac{1-\left(\frac{u}{4}\right)^{P}}{\frac{4}{u}-1}+\\
\left.+\frac{1}{G \, A^{*} \, v} \frac{1-\left(\frac{v}{4}\right)^{P}}{\frac{4}{v}-1}+\frac{2}{G \, A^{*} \, v} C_{P, \frac{v}{4}}\left(\sigma\left(z^{*}(w)\right)\right)\right)
\end{gathered}
\end{equation}

\item For $a\neq16$, $a\neq8$, $u\neq4$ and $v=4$ (terms (\ref{eq:43a}), (\ref{eq:46a}) and (\ref{eq:48b})):
\begin{equation} \label{eq:50}
\begin{gathered}
H \delta_{P, w}^{R}=2 \, c \, s \, L\left(\frac{10 \, L^{2}}{3 \, E \, I} \, \frac{1-\left(\frac{a}{16}\right)^{P}}{16-a}-\frac{20 \, L^{2}}{3 \, E \, I \, a} C_{P, \frac{a}{16}}\left(\sigma\left(z^{*}(w)\right)\right)-2 \, L^{2} \frac{\left(\frac{1}{2}\right)^P-\left(\frac{a}{16}\right)^{P}}{E \, I(8-a)}+\right.\\
+\frac{L^{2}}{2^{P-2} \, E \, I \, a} C_{P, \frac{a}{8}}\left(\sigma\left(z^{*}(w)\right)\right)+\frac{2}{E \, A \, u} C_{P, \frac{u}{4}}\left(\sigma\left(z^{*}(w)\right)\right)-\frac{1}{E \, A \, u} \frac{1-\left(\frac{u}{4}\right)^{P}}{\frac{4}{u}-1}+\\
\left.+\frac{P}{4 \, G \, A^{*}}-\frac{1}{2 \, G \, A^{*} } C_{P, 1}\left(\sigma\left(z^{*}(w)\right)\right)\right)
\end{gathered}
\end{equation}

\item For $a\neq16$, $a\neq8$, $u=4$ and $v\neq4$ (terms (\ref{eq:43a}), (\ref{eq:46b}) and (\ref{eq:48a})):
\begin{equation} \label{eq:51}
\begin{gathered}
H \delta_{P, w}^{R}=2 \, c \, s \, L\left(\frac{10 \, L^{2}}{3 \, E \, I} \, \frac{1-\left(\frac{a}{16}\right)^{P}}{16-a}-\frac{20 \, L^{2}}{3 \, E \, I \, a} C_{P, \frac{a}{16}}\left(\sigma\left(z^{*}(w)\right)\right)-2 \, L^{2} \frac{\left(\frac{1}{2}\right)^P-\left(\frac{a}{16}\right)^{P}}{E \, I(8-a)}+\right.\\
+\frac{L^{2}}{2^{P-2} \, E \, I \, a} C_{P, \frac{a}{8}}\left(\sigma\left(z^{*}(w)\right)\right)+\frac{1}{2 \, E \, A} C_{P, 1}\left(\sigma\left(z^{*}(w)\right)\right)-\frac{P}{4 \, E \, A}+\\
\left.+\frac{1}{G \, A^{*} \, v} \frac{1-\left(\frac{v}{4}\right)^{P}}{\frac{4}{v}-1}+\frac{2}{G \, A^{*} \, v} C_{P, \frac{v}{4}}\left(\sigma\left(z^{*}(w)\right)\right)\right)
\end{gathered}
\end{equation}

\item For $a\neq16$, $a\neq8$, $u=4$ and $v=4$ (terms (\ref{eq:43a}), (\ref{eq:46b}) and (\ref{eq:48b})):
\begin{equation} \label{eq:52}
\begin{gathered}
H\delta_{P, w}^{R}=2 \, c \, s \, L\left(\frac{10 \, L^{2}}{3 \, E \, I} \, \frac{1-\left(\frac{a}{16}\right)^{P}}{16-a}-\frac{20 \, L^{2}}{3 \, E \, I \, a} C_{P, \frac{a}{16}}\left(\sigma\left(z^{*}(w)\right)\right)-2 \, L^{2} \frac{\left(\frac{1}{2}\right)^P-\left(\frac{a}{16}\right)^{P}}{E \, I(8-a)}+\right.\\
+\frac{L^{2}}{2^{P-2} \, E \, I \, a} C_{P, \frac{a}{8}}\left(\sigma\left(z^{*}(w)\right)\right)+\frac{1}{2 \, E \, A} C_{P, 1}\left(\sigma\left(z^{*}(w)\right)\right)+\frac{P}{4 \, E \, A}+\\
\left.+\frac{P}{4 \, G \, A^{*}}-\frac{1}{2 \, G \, A^{*} } C_{P, 1}\left(\sigma\left(z^{*}(w)\right)\right)\right)
\end{gathered}
\end{equation}

\item For $a=8$, $u\neq4$ and $v\neq4$ (terms (\ref{eq:43b}), (\ref{eq:46a}) and (\ref{eq:48a})):
\begin{equation} \label{eq:53}
\begin{gathered}
H \delta_{P, w}^{R}=2 \, c \, s \, L\left(\frac{10 \, L^{2}}{3 \, E \, I} \, \frac{1-\left(\frac{1}{2}\right)^P}{8}-\frac{5 \, L^{2}}{6 \, E \, I} C_{P, \frac{1}{2}}\left(\sigma\left(z^{*}(w)\right)\right)-\frac{L^2 \, P}{2^{P+2}\, E \, I}+\frac{L^{2}}{2^{P+1} \, E \, I} C_{P, 1}\left(\sigma\left(z^{*}(w)\right)\right)+\right.\\
\left.+\frac{2}{E \, A \, u} C_{P, \frac{u}{4}}\left(\sigma\left(z^{*}(w)\right)\right)-\frac{1}{E \, A \, u} \frac{1-\left(\frac{u}{4}\right)^{P}}{\frac{4}{u}-1}+\frac{1}{G \, A^{*} \, v} \frac{1-\left(\frac{v}{4}\right)^{P}}{\frac{4}{v}-1}+\frac{2}{G \, A^{*} \, v} C_{P, \frac{v}{4}}\left(\sigma\left(z^{*}(w)\right)\right)\right)
\end{gathered}
\end{equation}

\item For $a=8$, $u\neq4$ and $v=4$ (terms (\ref{eq:43b}), (\ref{eq:46a}) and (\ref{eq:48b})):
\begin{equation} \label{eq:54}
\begin{gathered}
H \delta_{P, w}^{R}=2 \, c \, s \, L\left(\frac{10 \, L^{2}}{3 \, E \, I} \, \frac{1-\left(\frac{1}{2}\right)^P}{8}-\frac{5 \, L^{2}}{6 \, E \, I} C_{P, \frac{1}{2}}\left(\sigma\left(z^{*}(w)\right)\right)-\frac{L^2 \, P}{2^{P+2}\, E \, I}+\frac{L^{2}}{2^{P+1} \, E \, I} C_{P, 1}\left(\sigma\left(z^{*}(w)\right)\right)+\right.\\
\left.+\frac{2}{E \, A \, u} C_{P, \frac{u}{4}}\left(\sigma\left(z^{*}(w)\right)\right)-\frac{1}{E \, A \, u} \frac{1-\left(\frac{u}{4}\right)^{P}}{\frac{4}{u}-1}+\frac{P}{4 \, G \, A^{*}}+\frac{1}{2 \, G \, A^{*} } C_{P, 1}\left(\sigma\left(z^{*}(w)\right)\right)\right)
\end{gathered}
\end{equation}

\item For $a=8$, $u=4$ and $v\neq4$ (terms (\ref{eq:43b}), (\ref{eq:46b}) and (\ref{eq:48a})):
\begin{equation} \label{eq:55}
\begin{gathered}
H \delta_{P, w}^{R}=2 \, c \, s \, L\left(\frac{10 \, L^{2}}{3 \, E \, I} \, \frac{1-\left(\frac{1}{2}\right)^P}{8}-\frac{5 \, L^{2}}{6 \, E \, I} C_{P, \frac{1}{2}}\left(\sigma\left(z^{*}(w)\right)\right)-\frac{L^2 \, P}{2^{P+2}\, E \, I}+\frac{L^{2}}{2^{P+1} \, E \, I} C_{P, 1}\left(\sigma\left(z^{*}(w)\right)\right)+\right.\\
\left.+\frac{1}{2 \, E \, A} C_{P, 1}\left(\sigma\left(z^{*}(w)\right)\right)-\frac{P}{4 \, E \, A}+\frac{1}{G \, A^{*} \, v} \frac{1-\left(\frac{v}{4}\right)^{P}}{\frac{4}{v}-1}-\frac{2}{G \, A^{*} \, v} C_{P, \frac{v}{4}}\left(\sigma\left(z^{*}(w)\right)\right)\right)
\end{gathered}
\end{equation}

\item For $a=8$, $u=4$ and $v=4$ (terms (\ref{eq:43b}), (\ref{eq:46b}) and (\ref{eq:48b})):
\begin{equation} \label{eq:56}
\begin{gathered}
H \delta_{P, w}^{R}=2 \, c \, s \, L\left(\frac{10 \, L^{2}}{3 \, E \, I} \, \frac{1-\left(\frac{1}{2}\right)^P}{8}-\frac{5 \, L^{2}}{6 \, E \, I} C_{P, \frac{1}{2}}\left(\sigma\left(z^{*}(w)\right)\right)-\frac{L^2 \, P}{2^{P+2}\, E \, I}+\frac{L^{2}}{2^{P+1} \, E \, I} C_{P, 1}\left(\sigma\left(z^{*}(w)\right)\right)+\right.\\
\left.+\frac{1}{2 \, E \, A} C_{P, 1}\left(\sigma\left(z^{*}(w)\right)\right)-\frac{P}{4 \, E \, A}+\frac{P}{4 \, G \, A^{*}}+\frac{1}{2 \, G \, A^{*} } C_{P, 1}\left(\sigma\left(z^{*}(w)\right)\right)\right)
\end{gathered}
\end{equation}

\item For $a=16$, $u\neq4$ and $v\neq4$ (terms (\ref{eq:43c}), (\ref{eq:46a}) and (\ref{eq:48a})):
\begin{equation} \label{eq:57}
\begin{gathered}
H \delta_{P, w}^{R}=2 \, c \, s \, L\left(\frac{5 \, P \, L^{2}}{24 \, E \, I}-\frac{5 \, L^{2}}{12 \, E \, I} C_{P, 1}\left(\sigma\left(z^{*}(w)\right)\right)+\frac{L^2}{E\, I}\frac{\left(\frac{1}{2}\right)^P-1}{4}+\frac{L^{2}}{2^{P+2} \, E \, I} C_{P, 2}\left(\sigma\left(z^{*}(w)\right)\right)+\right.\\
\left.+\frac{2}{E \, A \, u} C_{P, \frac{u}{4}}\left(\sigma\left(z^{*}(w)\right)\right)-\frac{1}{E \, A \, u} \frac{1-\left(\frac{u}{4}\right)^{P}}{\frac{4}{u}-1}+\frac{1}{G \, A^{*} \, v} \frac{1-\left(\frac{v}{4}\right)^{P}}{\frac{4}{v}-1}+\frac{2}{G \, A^{*} \, v} C_{P, \frac{v}{4}}\left(\sigma\left(z^{*}(w)\right)\right)\right)
\end{gathered}
\end{equation}

\item For $a=16$, $u\neq4$ and $v=4$ (terms (\ref{eq:43c}), (\ref{eq:46a}) and (\ref{eq:48b})):
\begin{equation} \label{eq:58}
\begin{gathered}
H \delta_{P, w}^{R}=2 \, c \, s \, L\left(\frac{5 \, P \, L^{2}}{24 \, E \, I}-\frac{5 \, L^{2}}{12 \, E \, I} C_{P, 1}\left(\sigma\left(z^{*}(w)\right)\right)+\frac{L^2}{E\, I}\frac{\left(\frac{1}{2}\right)^P-1}{4}+\frac{L^{2}}{2^{P+2} \, E \, I} C_{P, 2}\left(\sigma\left(z^{*}(w)\right)\right)+\right.\\
+\left.\frac{2}{E \, A \, u} C_{P, \frac{u}{4}}\left(\sigma\left(z^{*}(w)\right)\right)-\frac{1}{E \, A \, u} \frac{1-\left(\frac{u}{4}\right)^{P}}{\frac{4}{u}-1}+\frac{P}{4 \, G \, A^{*}}-\frac{1}{2 \, G \, A^{*} } C_{P, 1}\left(\sigma\left(z^{*}(w)\right)\right)\right)
\end{gathered}
\end{equation}

\item For $a=16$, $u=4$ and $v\neq4$ (terms (\ref{eq:43c}), (\ref{eq:46b}) and (\ref{eq:48a})):
\begin{equation} \label{eq:59}
\begin{gathered}
H \delta_{P, w}^{R}=2 \, c \, s \, L\left(\frac{5 \, P \, L^{2}}{24 \, E \, I}-\frac{5 \, L^{2}}{12 \, E \, I} C_{P, 1}\left(\sigma\left(z^{*}(w)\right)\right)+\frac{L^2}{E\, I}\frac{\left(\frac{1}{2}\right)^P-1}{4}+\frac{L^{2}}{2^{P+2} \, E \, I} C_{P, 2}\left(\sigma\left(z^{*}(w)\right)\right)+\right.\\
\left.+\frac{1}{2 \, E \, A} C_{P, 1}\left(\sigma\left(z^{*}(w)\right)\right)-\frac{P}{4 \, E \, A}+\frac{1}{G \, A^{*} \, v} \frac{1-\left(\frac{v}{4}\right)^{P}}{\frac{4}{v}-1}-\frac{2}{G \, A^{*} \, v} C_{P, \frac{v}{4}}\left(\sigma\left(z^{*}(w)\right)\right)\right)
\end{gathered}
\end{equation}

\item For $a=16$, $u=4$ and $v=4$ (terms (\ref{eq:43c}), (\ref{eq:46b}) and (\ref{eq:48b})):
\begin{equation} \label{eq:60}
\begin{gathered}
H \delta_{P, w}^{R}=2 \, c \, s \, L\left(\frac{5 \, P \, L^{2}}{24 \, E \, I}-\frac{5 \, L^{2}}{12 \, E \, I} C_{P, 1}\left(\sigma\left(z^{*}(w)\right)\right)+\frac{L^2}{E\, I}\frac{\left(\frac{1}{2}\right)^P-1}{4}+\frac{L^{2}}{2^{P+2} \, E \, I} C_{P, 2}\left(\sigma\left(z^{*}(w)\right)\right)+\right.\\
\left.+\frac{1}{2 \, E \, A} C_{P, 1}\left(\sigma\left(z^{*}(w)\right)\right)-\frac{P}{4 \, E \, A}+\frac{P}{4 \, G \, A^{*}}+\frac{1}{2 \, G \, A^{*} } C_{P, 1}\left(\sigma\left(z^{*}(w)\right)\right)\right)
\end{gathered}
\end{equation}
\end{itemize}

\subsection{\bf{Horizontal Displacements in the Structure with Infinite Levels}} \label{sec:hinfty}

Once again, the limit $P\rightarrow\infty$ is taken to study the structure with infinite levels, where $P$ is the number of levels of the finite structure. The horizontal displacement $H\delta_{P,w}^R$ of the end node $w$ is obtained by adding up the expressions (\ref{eq:42}), (\ref{eq:45}) and (\ref{eq:47}). To facilitate the calculations when taking $P\rightarrow\infty$, it is appropriate to express $H\delta_{P,w}^R$ as:
\begin{equation}
\begin{gathered} \label{eq:61}
H \delta_{P, w}^{R}=2 \, c \, s \, L \sum_{i=1}^{P}\left(1-2 \rho_{i}\left(\sigma\left(z^{*}(w)\right)\right)\right)\left(\frac{10 \, L^{2}}{3 \, E \, I \, a}\left(\frac{a}{16}\right)^{i}+\right.\\
-\frac{L^{2}}{2^{P-1} \, E \, I \, a}\left(\frac{a}{8}\right)^{i}-\frac{1}{E \, A \, u}\left(\frac{u}{4}\right)^{i}+\left.+\frac{1}{G \, A^{*} \, v}\left(\frac{v}{4}\right)^{i}\right)
\end{gathered}
\end{equation}

In the second parenthesis of the expression (\ref{eq:61}), there are positive and negative terms that can cancel each other out as shown below:

\begin{itemize}
\item The terms $\frac{L^{2}}{2^{P-1} \, E \, I \, a}\left(\frac{a}{8}\right)^{i}$ and $\frac{1}{G \, A^{*} \, v}\left(\frac{v}{4}\right)^{i}$ cancel each other out for $a = 2 \, v, A^{*}=\frac{2^{P} \, E \, I}{G \, L^{2}}$, so there is not functional dependence on $v$. When taking $P\rightarrow\infty$ it results $A^\ast\rightarrow\infty$ (infinite shear area), which is absurd.
\item The terms $\frac{10 \, L^{2}}{3 \, E \, I \, a}\left(\frac{a}{16}\right)^{i}$ and $\frac{1}{E \, A \, u}\left(\frac{u}{4}\right)^{i}$ cancel each other out for $a=4\, u, A=\frac{6\, I}{5\, L^2}$, so there is not functional dependence on $u$. Therefore, in this case, when taking $P\rightarrow\infty$, it results, according to Eqs. (\ref{eq:49}) - (\ref{eq:60}):
\begin{subequations}
\begin{empheq}[left=\empheqlbrace]{align}
\begin{aligned}[t]  \label{eq:62a}
	    & H\delta_{P\rightarrow\infty,w}^R=\frac{2 \, c \, s \, L}{G \, A^*}\left(\frac{1}{4-v}-\frac{2}{v}C_\frac{v}{4}\left(\sigma\left(z^{*}(w)\right)\right)\right) \text{if} \quad 1< a<16, v<4
           \end{aligned}\\ 
\begin{aligned}[t]  \label{eq:62b}
	    & H\delta_{P\rightarrow\infty,w}^R \quad \text{diverges if} \quad a\geq16 \text{ or } v\geq4
           \end{aligned}
\end{empheq}
\end{subequations}

\item The terms $\frac{1}{E \, A\, u}\left(\frac{u}{4}\right)^i$ and $\frac{1}{G \, A^*\, v}\left(\frac{v}{4}\right)^i$ cancel each other out for $u=v, E \, A \, u=G \, A^* \, v$, so there is not functional dependence on $u$ and $v$. Therefore, when taking $P\rightarrow\infty$ it results, according to Eqs. (\ref{eq:49}) - (\ref{eq:60}):
\begin{subequations}
\begin{empheq}[left=\empheqlbrace]{align}
\begin{aligned}[t]  \label{eq:63a}
	    & H\delta_{P\rightarrow\infty,w}^R=\frac{20 \, c \, s \, L^3}{3 \, E \, I}\left(\frac{1}{16-a}-\frac{2}{a}C_\frac{a}{16}\left(\sigma\left(z^{*}(w)\right)\right)\right) \text{if} \quad 1<a<8, 8<a<16
           \end{aligned}\\ 
\begin{aligned}[t]  \label{eq:63b}
	    & H\delta_{P\rightarrow\infty,w}^R=\frac{5 \, c \, s \, L^3}{3 \, E \, I}\left(\frac{1}{2}-C_\frac{1}{2}\left(\sigma\left(z^{*}(w)\right)\right)\right) \text{if} \quad a=8
           \end{aligned}\\ 
\begin{aligned}[t]  \label{eq:63c}
	    & H\delta_{P\rightarrow\infty,w}^R \quad \text{diverges if} \quad a\geq16
           \end{aligned}
\end{empheq}
\end{subequations}
\end{itemize}

If there are no cancellations between the terms in the second parenthesis of Eq. (\ref{eq:61}), when taking $P\rightarrow\infty$ it results, according to Eqs. (\ref{eq:49}) - (\ref{eq:60}): 

\begin{subequations}
\begin{empheq}[left=\empheqlbrace]{align}
\begin{aligned}[t] \label{eq:64a}
	    & H\delta_{P\rightarrow\infty,w}^R=2 \, c \, s \, L \left(\frac{10 \, L^2}{3 \, E \, I(16-a)}-\frac{20 \, L^2}{3 \, E \, I \, a}C_\frac{a}{16}\left(\sigma\left(z^{*}(w)\right)\right)+\frac{2}{E \, A \, u}C_\frac{u}{4}\left(\sigma\left(z^{*}(w)\right)\right)+\right.\\
	    & \left.-\frac{1}{E\, A(4-u)}+\frac{1}{G\, A^*(4-v)}-\frac{2}{G\, A^* \, v}C_\frac{v}{4}\left(\sigma\left(z^{*}(w)\right)\right)\right) \text{if} \quad 1<a<16, u<4, v<4 
           \end{aligned}\\ 
\begin{aligned}[t] \label{eq:64b}
	    & H\delta_{P\rightarrow\infty,w}^R \quad \text{diverges if} \quad a\geq16, u\geq4 \text{ and/or } v\geq4 
           \end{aligned}
\end{empheq}
\end{subequations}

Where
\begin{equation} \label{eq:65}
C_t(x) =\lim_{P\rightarrow\infty} C_{P,t}(x) = \lim_{P\rightarrow\infty} \sum_{k=1}^{P} \rho_k(x)\, t^k =\sum_{k=1}^{\infty} \rho_k(x)\, t^k
\end{equation}

\subsection{\bf{Horizontal Displacements and the $\beta$-Cantor Function}}

Let $x$ be a element of the $\beta$-Cantor set $\Omega_\beta$, such that
\begin{equation}
\Omega_\beta=\left\{\frac{1+\beta}{1-\beta}\sum_{k=1}^{\infty} \alpha_k \left(\frac{1-\beta}{2}\right)^k \bigg| \alpha_k\in \{0,1\}, k=0,1,\ldots \right\} \nonumber,
\end{equation}
the $\beta$-Cantor function $f_\beta$ holds \cite{21}:
\begin{equation} \label{eq:66}
f_\beta(x)=f_\beta\left(\frac{1+\beta}{1-\beta}\sum_{k=1}^{\infty} \alpha_k \left(\frac{1-\beta}{2}\right)^k \right)=\sum_{k=1}^{\infty}\frac{\alpha_k}{2^k} 
\end{equation}

Therefore:
\begin{equation} \label{eq:67}
f_\beta^{-1}\left(\sum_{k=1}^{\infty}\frac{\alpha_k}{2^k}\right)= \frac{1+\beta}{1-\beta}\sum_{k=1}^{\infty} \alpha_k \left(\frac{1-\beta}{2}\right)^k
\end{equation}

\begin{figure}[h]
\centering
\includegraphics[width=10cm]{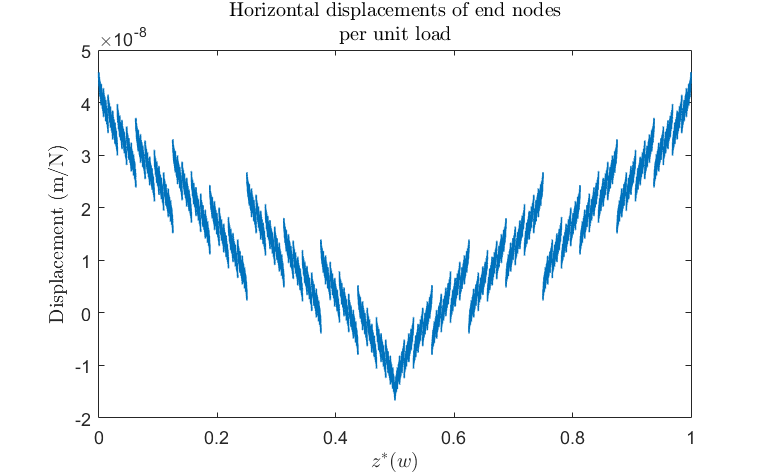}
\caption{Combination of inverses of $\beta$-Cantor functions for the horizontal displacements of the end nodes per unit load versus the end nodes positions $z^*(w)$ for a structure with infinite levels. A positive value indicates leftward displacement for $z^*(w)<0.5$ and rightward for $z^*(w)\geq0.5$. The horizontal displacements are plotted for $\theta=60^{\circ}$, $E=10^{10}$ N/m\textsuperscript{2}, $G=5\cdot10^8$ N/m\textsuperscript{2}, $L=0.5$ m, $I=3.1416\cdot10^{-4}$ m\textsuperscript{4}, $A=3.1416\cdot10^{-2}$ m\textsuperscript{2}, $A^*=2.8274\cdot10^{-2}$ m\textsuperscript{2}, $a=9$, $u=3$ and $v=3$.}
\label{fig:hdesp}
\end{figure}

Let $w$ be the end node located at $z^*(w)$ such that
\begin{equation}
\sigma\left(z^*(w)\right)=\sum_{k=1}^{\infty}\frac{\rho_k\left(\sigma\left(z^*(w)\right)\right)}{2^k} \nonumber
\end{equation}
as expressed in Sec. \ref{sec:hfinite}. Considering Eq. (\ref{eq:66}), it results:
\begin{equation}
\begin{gathered}
f_\beta^{-1}\left(\sigma\left(z^*(w)\right)\right)=f_\beta^{-1}\left(\frac{\rho_k\left(\sigma\left(z^*(w)\right)\right)}{2^k}\right)=
\frac{1+\beta}{1-\beta}\sum_{k=1}^{\infty} \rho_k\left(\sigma\left(z^*(w)\right)\right) \left(\frac{1-\beta}{2}\right)^k \nonumber
=\frac{1-t}{t}\sum_{k=1}^{\infty}\rho_k\left(\sigma\left(z^*(w)\right)\right)\, t^k\Rightarrow\\
\Rightarrow f_\beta^{-1}\left(\sigma\left(z^*(w)\right)\right)=\frac{1-t}{t}C_t\left(\sigma\left(z^*(w)\right)\right)
\end{gathered}
\end{equation}
with $t=\frac{1-\beta}{2}$ and $C_t$ given by Eq. (\ref{eq:65}). Therefore, for any end node $w$, the function $\frac{1-t}{t}C_t\left(\sigma\left(z^*(w)\right)\right)$ showed in Sec. \ref{sec:hinfty} is the inverse of the $\beta$-Cantor function $f_\beta$. The interpretation of this result is that the horizontal displacements in the infinite structure $H\delta_{P\rightarrow\infty,w}^R$ given by Eqs. (\ref{eq:62a}), (\ref{eq:63a}), (\ref{eq:63b}) and (\ref{eq:64a}) are the linear combination of three inverse $\beta$-Cantor functions (see Fig. \ref{fig:hdesp}). Note that the Devil's Staircase associated to the classic Cantor set is a particular case of the $\beta$-Cantor function for $\beta=1/3$. 
\\

%% file: iteraciones.tex
\section{\bf{DISPLACEMENTS in all NODES: ITERATIONS LEADING to the FRACTALS}}\label{sec:iterations}

\begin{figure}[h]
     \centering
     \begin{subfigure}
         \centering
         \includegraphics[width=6cm]{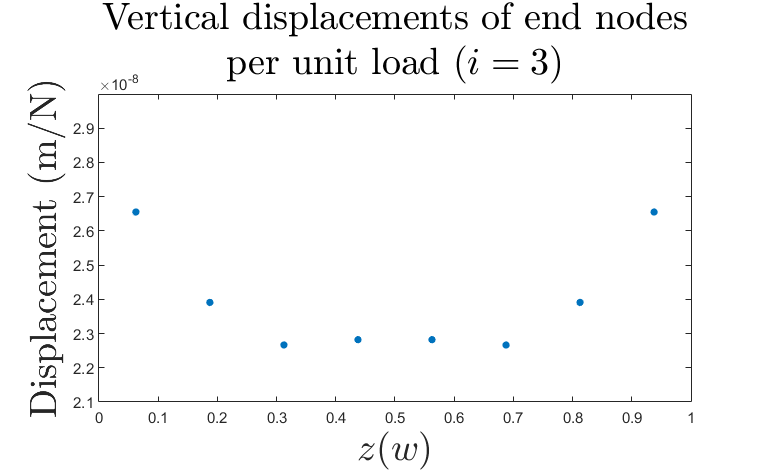}
     \end{subfigure}
     \hspace{-0.5cm}
     \begin{subfigure}
         \centering
         \includegraphics[width=6cm]{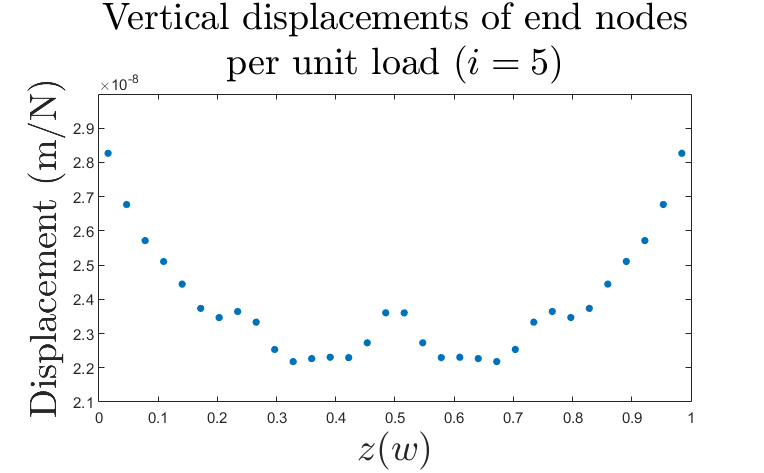}\\
     \end{subfigure}
     \par\bigskip
     \begin{subfigure}
         \centering
         \includegraphics[width=6cm]{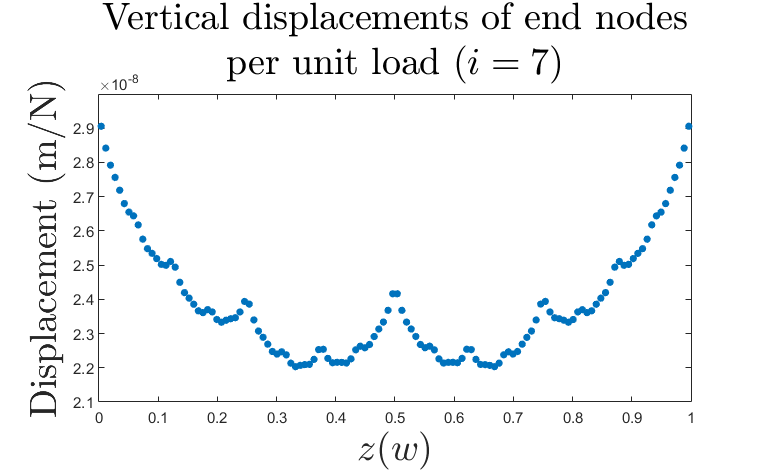}
     \end{subfigure}
     \hspace{-0.5cm}
     \begin{subfigure}
         \centering
         \includegraphics[width=6cm]{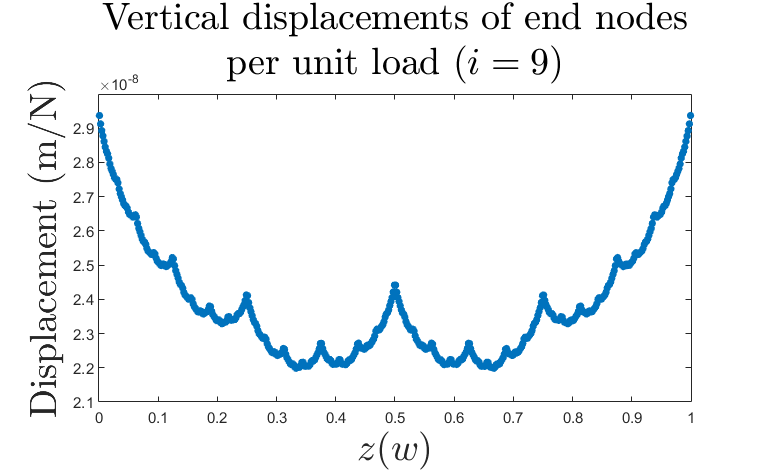}
     \end{subfigure}
     \par
        \caption{Vertical displacements per unit load versus the end nodes positions $z(w)$ in a binary tree structure. The higher the level $i$, the better the approximation to the Takagi curve.}
        \label{fig:vdespinc}
\end{figure}

\begin{figure}[H]
     \centering
     \begin{subfigure}
         \centering
         \includegraphics[width=6cm]{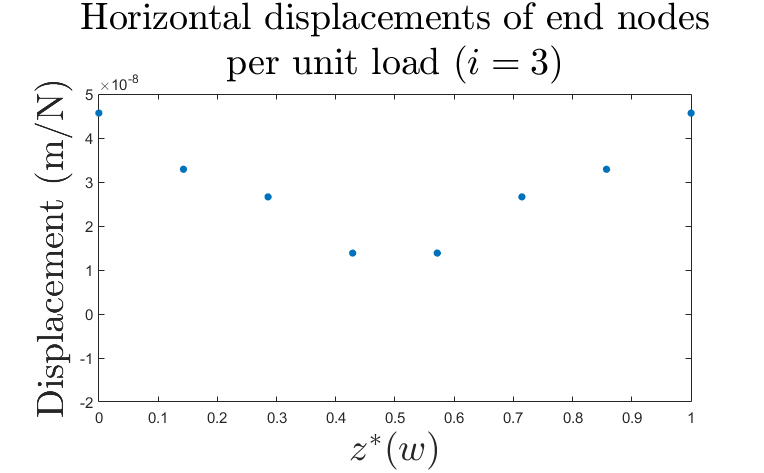}
     \end{subfigure}
     \hspace{-0.5cm}
     \begin{subfigure}
         \centering
         \includegraphics[width=6cm]{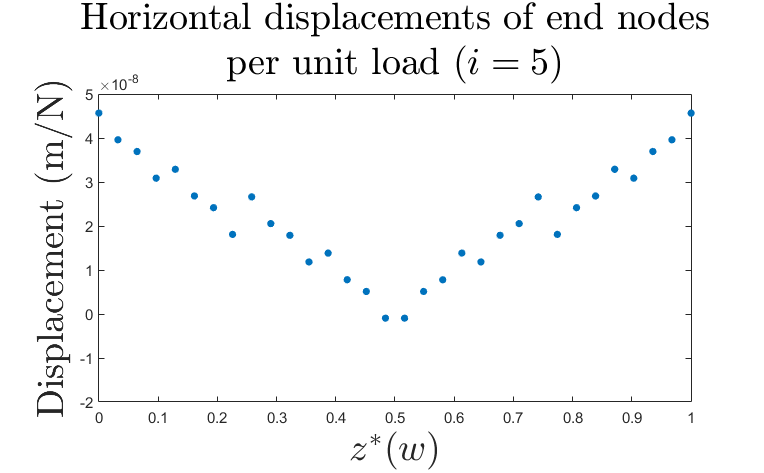}\\
     \end{subfigure}
     \par\bigskip
     \begin{subfigure}
         \centering
         \includegraphics[width=6cm]{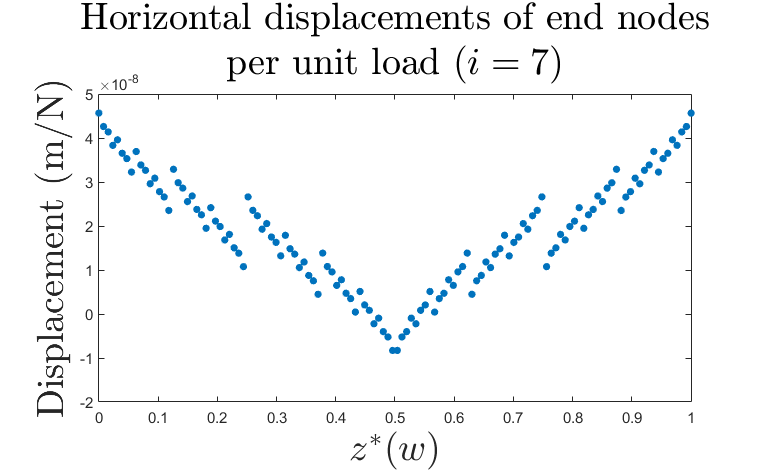}
     \end{subfigure}
     \hspace{-0.5cm}
     \begin{subfigure}
         \centering
         \includegraphics[width=6cm]{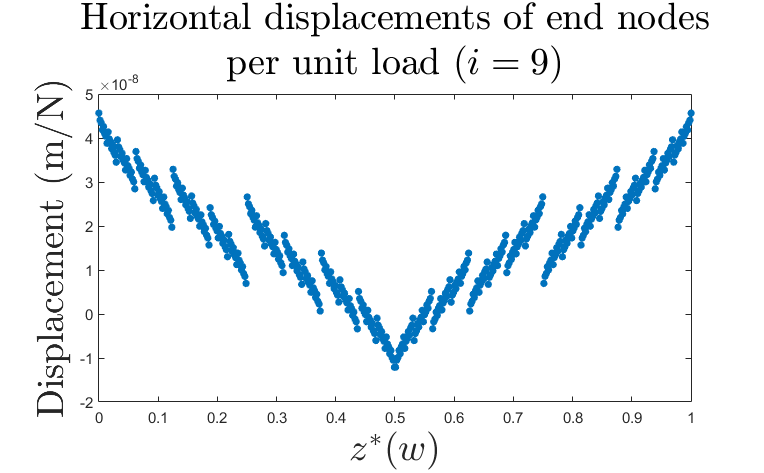}
     \end{subfigure}
     \par
        \caption{Horizontal displacements per unit load versus the end nodes positions $z^*(w)$ in a binary tree structure. The higher the level $i$, the better the approximation to the linear combination of inverses of $\beta$-Cantor functions.}
        \label{fig:hdespinc}
\end{figure}

Let us consider the structure defined in Sec. \ref{sec:model}. In this structure, the nodes of a level $i$, $i\leq P$, receive a load of $1/2^i$ each, so they behave like the end nodes of a structure of $i$ levels that receives a total load of value 1. If $n$ is a node from level $i$, its total vertical displacement $V\delta_{i,n}^R$ is given by Eqs. (\ref{eq:24})-(\ref{eq:31}) and its total horizontal displacement $H\delta_{i,n}^R$ by Eqs. (\ref{eq:49})-(\ref{eq:60}), taking $w=n$ and $P=i$ in all of them. These equations represent the iterations that lead to the Takagi curve and the linear combination of inverses of $\beta$-Cantor functions, as can be seen in Figs. \ref{fig:vdespinc} and \ref{fig:hdespinc}, respectively.
\\

The higher the level $i$, the better the iterations approximate their corresponding fractals. Therefore, the nodes from a level $i$ of a binary tree structure deform vertically following the $i$-th iteration of a Takagi curve and horizontally following the $i$-th iteration of a linear combination of inverses of $\beta$-Cantor functions.

%% file: dimension.tex
\section{\bf{The RELATIONSHIP Between TAKAGI and the $\beta$-CANTOR FUNCTION via the FRACTAL DIMENSION}}\label{sec:dimensions}
\
The fractal dimension $D_\Psi$ of a Takagi curve $\Psi_\frac{a}{16}$, which appears in Eq. (\ref{eq:33}) due to the bending moments, is given by \cite{22}:
\begin{equation} \label{eq:68}
 D_\Psi = \frac{\log(a/4)}{\log(2)}
\end{equation}

On the other hand, the fractal dimension $D_f$ of the $\beta$-Cantor function $f_\beta$ showed in Eq. (\ref{eq:66}) is given by \cite{21}:
\begin{equation}
D_f = \frac{-\log(2)}{\log\left(\frac{1-\beta}{2}\right)} \nonumber
\end{equation}

Therefore, the fractal dimension $D_C$ of the inverse of a $\beta$-Cantor function $C_t$ with $t=\frac{1-\beta}{2}$ is:
\begin{equation} \label{eq:69}
D_C = \frac{-\log(2)}{\log(t)}
\end{equation}

Thus, according to Eq. (\ref{eq:69}), the fractal dimension of the function $C_\frac{a}{16}$ that appears in Eqs. (\ref{eq:63a}) and (\ref{eq:64a}) due to the bending moments is: 
\begin{equation} \label{eq:70}
D_C = \frac{-\log(2)}{\log(a/16)}
\end{equation}

Finally, according to Eqs. (\ref{eq:68}) and (\ref{eq:70}), the fractal dimensions $D_\Psi$ and $D_C$ of the Takagi curve and the inverse of the $\beta$-Cantor function that appears in the vertical displacements (Eq. (\ref{eq:33})) and horizontal displacements (Eqs. (\ref{eq:63a}) and (\ref{eq:64a})) due to the bending moments, are related by:
\begin{equation}
D_\Psi+\frac{1}{D_C}=2
\end{equation}

%% file: conclusiones.tex
\section{\bf{CONCLUSIONS}}\label{sec:conclusions}
\
Fractals have been derived by using natural laws and not algorithms.

We have studied the vertical and horizontal deformations of the nodes in a loaded binary tree by applying continuum elasticity theory. As a consequence of this theory, two fractals emerge: the first one is associated to the vertical displacements of the end nodes and the second one to the horizontal displacements. In the first case, the Takagi curve appears, and in the second case, it is a linear combination of inverses of $\beta$-Cantor functions. In addition, a link between the fractal dimension of both fractals has been found in the studied structure. This result closes the gap between the Platonic world of Mathematics and the real world and shows how fractals emerge from the laws of Nature, without being forced by any mathematical algorithm. 

The results are more profound than finding fractals from the laws of Mechanics. The first is that a tree can generate a fractal in its crown by varying the inertia of its branches without changing its geometry. The inertia reduction results in the fractal structure of both vertical (Tagaki curve) and horizontal ($\beta$-Cantor function) displacements when the number of levels of the structure tends to infinity.

Let us emphasize that that this type of ﬁnite binary tree structures are increasingly frequent in mega-structures, where the aim is to eliminate as much of the pillars as possible, as can be seen in new airport terminals and train stations, among others. The key point here is that the calculation of the displacements is given explicitly by the mentioned expressions without the need to use structure calculation programs.

The second and unexpected finding is that fractals have appeared in pairs in a binary tree: Takagi curve and $\beta$-Cantor function. It is interesting to note that the eigenvalue distribution of a binary tree of spring-connected masses exhibits a Devil’s Staircase self-similarity \cite{14} (the Devil's Staircase of the classic Cantor set is a particular case of the $\beta$-Cantor function for $\beta=1/3$). The presence of the $\beta$-Cantor function both in our work and in the mentioned reference \cite{14} makes us wonder if the $\beta$-Cantor function will go hand in hand with binary trees; and if the $\beta$-Cantor function appears, so will the Takagi curve. Note that the elastic deformation of the branches can be seen as a compressed spring, so we can conjecture that those equivalent cases to binary tree of spring-connected masses, or formally identical in their mathematical formulation, will show both Takagi curve and $\beta$-Cantor function. The presence of binary trees is abundant in the literature \cite{15,16,17,18,19}, so the phenomenon may be more common than the two cases presented here, specially taking into account, as mathematically proved, that fractal dimension of both fractals are related in the binary tree structure.